\documentclass[12pt,reqno,intlim]{amsart}

\usepackage{amssymb,amsmath}

\usepackage[arrow,matrix,curve]{xy}

\newdir{ >}{{}*!/-5pt/\dir{>}}

\newcommand{\nc}{\newcommand}

\nc{\bC}{\bold{C}} \nc{\bN}{\Bbb{N}} \nc{\cF}{\mathcal{F}}
\nc{\cE}{\mathcal{E}} \nc{\cR}{\mathcal{R}} \nc{\cM}{\mathcal{M}}
\nc{\al}{\alpha} \nc{\bt}{\beta} \nc{\gm}{\gamma} \nc{\dl}{\delta}
\nc{\om}{\omega} \nc{\sg}{\sigma} \nc{\Sg}{\Sigma} \nc{\vf}{\varphi}
\nc{\ve}{\varepsilon} \nc{\os}{\overset} \nc{\ol}{\overline}
\nc{\ul}{\underline} \nc{\us}{\underset} \nc{\sbs}{\subset}
\nc{\bsl}{\backslash} \nc{\Ra}{\Rightarrow}
\nc{\lra}{\longrightarrow} \nc{\all}{\allowdisplaybreaks}

\nc{\Codes}{\operatorname{{\bold{Codes}}}}
\nc{\RegMono}{\operatorname{\mathcal{R}{\rm{eg}\mathcal{M}{\rm{ono}\!}}}}
\nc{\RegEpi}{\operatorname{\mathcal{R}{\rm{eg}\mathcal{E}{\rm{pi}\!}}}}
\nc{\Mn}{\operatorname{\mathcal{M}{\rm{ono}\!}}}
\nc{\Ep}{\operatorname{\mathcal{E}{\rm{pi}\!}}}
\nc{\Rg}{\operatorname{\mathcal{R}{\rm{eg}\!}}}
\nc{\Ob}{\operatorname{Ob\!}}

\numberwithin{equation}{section}

\newtheorem{theo}{\ \ \ Theorem}[section]
\newtheorem{lem}[theo]{\ \ \ Lemma}
\newtheorem{prop}[theo]{\ \ \ Proposition}
\newtheorem{cor}[theo]{\ \ \ Corollary}

\theoremstyle{definition}

\theoremstyle{remark}
\newtheorem{rem}[theo]{\ \ \ Remark}

\begin{document}

\title[]
{On stable - projective and injective - costable decompositions of modules}

\author{Dali Zangurashvili}

\maketitle

\begin{abstract}
It is proved that, for a left hereditary ring, an arbitrary left module has a representation in the form of the direct sum of a stable left module and indecomposable projective left modules (if and only if an arbitrary left module has a representation in the form of the direct sum of a stable left module and a projective left module) if and only if the ring is left perfect and right coherent. In that case, the above-mentioned representations are unique up to isomorphism; the latter representation is also functorial.  The essential ingredient in the proofs of the above-mentioned statements is a certain purely categorical result. These statements, in particular, imply that, for any principal ideal domain that is not a field, the fundamental theorem on finitely generated modules over it can not be generalized to the case of all modules. Moreover, with the aid of the above-mentioned categorical approach, we give a new proof of the Zheng-Xu He's result asserting that any module of a ring has a unique up to isomorphism  injective-costable decomposition if and only if the ring is left hereditary and left Noetherian. The above-mentioned statements, in particular, imply that if the category of left modules over a left hereditary ring is Krull-Schmidt, then the ring is left Artinian.  Yet another criterion for a ring to be left hereditary, left perfect and right coherent (resp. left hereditary left Noetherian) found in the paper requires that the pair $(Stable$ $modules$, $Projective$ $modules)$ (resp. $(Injective$ $Modules, Costable$ $ Modules)$) of module classes be a pre-torsion theory. This implies that the pair $(Stable$ $modules$, $Projective$ $modules)$ is a torsion theory if and only if the ring is left hereditary and the injective envelope of the ring, viewed as a left module over itself, is projective. 

\vskip+2mm
\noindent{\bf Key words and phrases}: left hereditary ring, left perfect right coherent ring, left Noetherian ring, (co)stable module, projective module, injective module, Freyd adjoint functor theorem, radical, (pre-)torsion theory, Krull-Schmidt category.
\vskip+2mm

\end{abstract}

\section{Introduction}

Throughout the paper `ring' means an associative ring with unit, and `module' means a left module over the ring.

Recall that a module is called stable if it has no nonzero projective summands. If a ring is left hereditary, then the direct sum of an arbitrary small family of stable modules is stable (see Section 3 below). Moreover, each indecomposable module is obviously either stable or projective. Therefore, as it was observed by Alex Martsinkovsky, if a module has an indecomposable decomposition, then it can be represented as the direct sum of a stable module and indecomposable projective modules, and also as the direct sum of a stable module and a projective module (we call the corresponding decompositions stable-indecomposable projectives and stable-projective decompositions respectively). Similarly, if a module has a \textit{finite} indecomposable decomposition, then it can be represented as the direct sum of indecomposable injective modules and a costable module (i.e, a module that is i-reduced in the sense of \cite{Zh}), and also as the direct sum of an injective module and a costable module. The rings over which all modules have decompositions of the latter two kinds (i.e., the indecomposable injectives-costable and injective-costable decompositions) were characterized by He \cite{Zh}. Namely, He showed that these are precisely left Noetherian rings. Moreover, He proved that such decompositions exist and are unique up to isomorphism if and only if a ring is left hereditary and left Noetherian. 

In view of the above-mentioned facts, the problem of characterizing rings over which all modules have (unique up to isomorphism) stable-projective/stable-indecomposable projectives decompositions naturally arises. The aim of this paper is to give such a characterization in the case of left hereditary rings. We prove that, for such a ring, the following conditions are equivalent:\vskip+1mm
(i) \textit{any module can be represented as the direct sum of a stable module and indecomposable projective modules;}
\vskip+1mm
(ii) \textit{any module can be represented as the direct sum of a stable module and a projective module};
\vskip+1mm
(iii)  \textit{the ring is left perfect and right coherent.}
\vskip+1mm
We also show that when these conditions are satisfied,  \textit{both above-mentioned representations  are unique up to isomorphism. The representation in the condition (ii) is also functorial.}\vskip
+1mm

The essential ingredient of our proof is a purely categorical result -- the old theorem on the reflectivity of some full replete subcategories of a complete well-powered co-well-powered category. We give the generalized form of this theorem replacing the requirement of well-poweredness of the ground category by the requirement of the existence of a factorization system satisfying a certain natural condition, and, moreover, give a certain form of the converse statement. The above-mentioned characterization of left hereditary left perfect right coherent rings is obtained  by applying these statements, the criterion for the direct product of projective modules to be projective by Chase \cite{C}, as well as the result on indecomposable decompositions of projective modules by Anderson and Fuller \cite{AF}. 

Note that the fundamental theorem on finitely generated modules over a principal ideal domain leads to stable-projective representations of such modules. Therefore, our criterion for the existence of such decompositions implies:
\vskip+1mm
- \textit{for any principal ideal domain that is not a field, the fundamental theorem can not be generalized to the case of all modules over this domain}.\vskip+2mm

Further, applying the dual of the `old theorem' and its converse, the criterion for the direct sum of injective modules to be injective by Bass \cite{Ba1}, as well as the result on indecomposable decompositions of injective modules by Matlis \cite{Ma}, we give a new proof of the He's characterization of left hereditary left Noetherian rings. Moreover, we show that the injective-costable representation of modules over such rings is functorial.


The above-mentioned results lead to a partial answer to the general problem: 
\vskip+1mm
- \textit{if any module over a left hereditary ring has a finite indecomposable decomposition, then the ring is left Artinian}.\vskip+1mm

With this regard, it is appropriate to mention here the result by Anderson and Fuller asserting that if any module over a ring has a decomposition that complements direct summands (and hence is indecomposable), then the ring is left Artinian \cite{AF}.\vskip+1mm

In the present paper, in addition to the above-mentioned  criteria for a left hereditary ring to be left perfect and right coherent, we give also some other  characterizations of such rings, as well as of left Noetherian rings. One of these characterizations
asserts that a ring is left hereditary, left perfect and right coherent if and only if the pair  
\begin{equation}
(Stable\; Modules, Projective \; Modules)
\end{equation} 
\noindent of module classes is a pre-torsion theory. A similar result is proved for left hereditary Noetherian rings: a ring is of this kind if and only if the pair 
\begin{equation}
(Injective\; Modules, Costable\; Modules)
\end{equation} 
\noindent of module classes is a pre-torsion theory. To this end, a certain idempotent radical on the category of modules over an arbitrary left hereditary ring is introduced.

Finally note that pair (1.2) is a torsion theory if and only if the ring is semisimple. As to pair (1.1), it is a torsion theory if and only if the ring is left hereditary and its injective envelope, viewed as a left module over itself, is projective. This statement follows from the main result of the paper \cite{WMJ} by Wu, Mochizuki, and Jans and the above-mentioned result on pair (1.1). 

The author expresses her gratitude to Alex Martsinkovsky and Mamuka Jibladze for valuable talks on the subject of the paper. 

Financial support from  Shota Rustaveli  Georgian National Science Foundation
(Ref.: FR-18-10849) is gratefully acknowledged.

\section{Notation and blanket assumptions}

Let $\Lambda$ be a ring, and $\Lambda$-$Mod$ be the category of $\Lambda$-modules.

Throughout the paper, when no confusion might arise, we use the same symbol for a class of modules and the corresponding full subcategory of the category $\Lambda$-$Mod$. Moreover, the symbol $_{\Lambda}\Lambda$ denotes the ring $\Lambda$, viewed as a left module over itself.
Further, for a module class $\mathcal{C}$,  the symbol $\mathcal{C}^{\bot}$ denotes the class of all modules $N$ such that $Hom_{\Lambda}(M,N)=0$, for any module $M$ from $\mathcal{C}$; and  $\mathcal{C}^{\top}$ denotes the class of all modules $M$ such that $Hom_{\Lambda}(M,N)=0$ for any module $N$ from $\mathcal{C}$. 

For a category $\mathcal{C}$, the symbol $Epi$ (resp. $Mono$) denotes the class of epimorphisms (resp. monomorphisms), while $StrongEpi$ (resp. $StrongMono$) denotes the class of strong epimorphisms (resp. strong monomorphisms) of $\mathcal{C}$.

\section{The classes of stable and costable modules}
In \cite{WMJ}, Wu, Mochizuki, and Jans noted that the class $\Lambda$-$Mod_{ZD}$ of modules with the zero duals is closed under homomorphic images, direct sums, and extensions (but, in general, is not closed under submodules). Hence it determines some idempotent radical $\textbf{R}$ on the category $\Lambda$-$Mod$. Note that the class $(\Lambda$-$Mod_{ZD})^{\perp}$ contains the class of torsionless modules, but, in general, does not coincide with it, as it follows from Theorem  
(on page 8) of \cite{WMJ}.

 In Lemma 2.6 of \cite{MZ}, Martsinkovsky and the present author observed that any module with the zero dual is stable, and the converse holds if a ring is left hereditary. Therefore, we obtain the following statement.
 \begin{lem}
If the direct sum of an arbitrary small family of stable $\Lambda$-modules is stable, then all these modules are stable. The converse is true if $\Lambda$ is left hereditary. In that case, the class $\Lambda$-$Mod_{St}$ is closed under (direct sums), homomorphic images and extensions.
\end{lem}
 Let $\Lambda$-$Mod_{St}$ be the class/category of stable modules, and $\Lambda$-$Mod_{Pr}$ be the class/category of projective modules. In view of Lemma 2.6 of \cite{MZ}, we also provide the following easily verified fact.
\begin{lem} Let $M$ be a module. For the following conditions, we have 
\begin{center}
(i)$\Leftrightarrow$(ii)$\Leftarrow$(iii). 
\end{center}
If $\Lambda$ is left hereditary, then the conditions (i)-(iii) are equivalent, and hence $\Lambda$-$Mod_{St}$=$(\Lambda$-$Mod_{Pr})^{\top}$:\vskip+2mm

(i) the module $M$ is stable;\vskip+1mm

(ii) there are no epimorphisms $M\twoheadrightarrow P$ with nonzero projective module $P$;\vskip+1mm

(iii) there are no nonzero homomorphisms $M\rightarrow P$ with projective $P$.
\end{lem}
Lemma 3.2 implies 
\begin{lem}
For an arbitrary ring $\Lambda$, the class $\Lambda$-$Mod_{St}$ of stable modules is closed under homomorphic images.
\end{lem}

Taking the above-mentioned facts into account, we obtain
\begin{prop}
Let $\Lambda$ be left hereditary. Then\vskip+1mm

(a) for any module $M$, $\textbf{R}(M)$ is the largest stable submodule of $M$; \vskip+1mm

(b) the pair  $(\Lambda$-$Mod_{St},(\Lambda$-$Mod_{St})^{\perp})$ of module classes is a pre-torsion theory on the category $\Lambda$-$Mod$;\vskip+1mm

(c) the subcategory $\Lambda$-$Mod_{St}$ of $\Lambda$-$Mod$ is coreflective;\vskip+1mm

(d) the subcategory $(\Lambda$-$Mod_{St})^{\perp}$ of $\Lambda$-$Mod$ is reflective;\vskip+1mm

(e)  the class $(\Lambda$-$Mod_{St})^{\perp}$ coincides with that of all modules $M$ for which $\textbf{R}(M)=0$.\vskip+1mm
\end{prop}

Dually to the notion of a stable module one can introduce the notion of a costable module. It obviously is equivalent to that of an i-reduced module by He \cite{Zh} (that is a generalization of the notion of a reduced Abelian group by Kaplansky \cite{K}).
\vskip+1mm
Let $\Lambda$-$Mod_{Cost}$ be the class/category of such modules, and $\Lambda$-$Mod_{Inj}$ be the class/category of injective modules. We have
\begin{lem} Let $M$ be a module. For the following conditions, we have 
\begin{center}
(i)$\Leftrightarrow$(ii)$\Leftarrow$(iii). 
\end{center}
If $\Lambda$ is left hereditary, then the conditions (i)-(iii) are equivalent, and hence $\Lambda$-$Mod_{Cost}$=$(\Lambda$-$Mod_{Inj})^{\bot}$:\vskip+2mm

(i) the module $M$ is costable;\vskip+1mm

(ii) there are no monomorphisms $I\rightarrowtail M$ with nonzero injective module $I$;\vskip+1mm

(iii) there are no nonzero homomorphisms $I\rightarrow M$ with injective module $I$.
\end{lem}

 Lemma 3.5 implies
\begin{lem}
The class $\Lambda$-$Mod_{Cost}$ of costable modules is closed under submodules.
If the direct sum of an arbitrary family of stable $\Lambda$-modules is costable, then all these modules are costable. If $\Lambda$ is left hereditary, then the class $\Lambda$-$Mod_{Cost}$ is closed under direct products and extensions.
\end{lem}

From Lemma 3.6 we obtain that the class $\Lambda$-$Mod_{Cost}$ determines some idempotent radical $\textbf{R}'$ on $\Lambda$-$Mod$.
\begin{prop}
Let $\Lambda$ be a left hereditary ring. Then\vskip+1mm

(a) for any module $M$, the submodule $\textbf{R}'(M)$ is smallest among all submodules $N$ with costable $M/N$; 
\vskip+1mm

(b) the pair $((\Lambda$-$Mod_{Cost})^{\top}$, $\Lambda$-$Mod_{Cost})$ of module classes is a pre-torsion theory;\vskip+1mm

(c) the subcategory $\Lambda$-$Mod_{Cost}$ of the category $\Lambda$-$Mod$ is reflective;\vskip+1mm

(d) the subcategory $(\Lambda$-$Mod_{Cost})^{\top}$ of the category $\Lambda$-$Mod$ is coreflective;\vskip+1mm

(e) the class $(\Lambda$-$Mod_{Cost})^{\top}$ coincides with that of all modules $M$ for which $\textbf{R}'(M)=M$.
\end{prop}
\vskip+1mm

 \section{The `old theorem'}
Recall the following well-known fact: \textit{let $\mathcal{C}$ be a complete well-powered and co-well-powered category, and $\mathcal{X}$ be a full replete subcategory of $\mathcal{C}$. If $\mathcal{X}$ is closed under products and subobjects, then it is reflective}. In \cite{Fr}, Freyd refers to this statement as the `old theorem', and, moreover, observes that it is an immediate consequence of his adjoint functor theorem. 

Below we give it in the generalized form where the requirement of well-poweredness of the category $\mathcal{C}$ is replaced by the requirement of the existence of a factorization system $(\mathbb{E},\mathbb{M})$ with the $\mathbb{E}\subseteq Epi$ on $\mathcal{C}$. We give also a certain form of the converse statement. These statements follow from Remark 3.1(ii) of \cite{Z}, but for the readers convenience, we give their direct proofs. Note that, in the particular case, where $\mathcal{C}$ is the category of topological spaces with the relevant classes $\mathbb{E}$ and $\mathbb{M}$, the obtained criterion is given in the paper \cite{Ke} by Kennison (see Theorem C). 

For the definition of a factorization system we refer the reader to, e.g., \cite{B}. Recall that if $(\mathbb{E},\mathbb{M})$ is a factorization system on $\mathcal{C}$, then the intersection $\mathbb{E}\cap\mathbb{M}$ coincides with the class of all isomorphisms of $\mathcal{C}$. Moreover, if $\mathbb{E}\subseteq Epi$, then the morphism class $\mathbb{M}$ satisfies the left cancellation property, i.e. if $gf\in \mathbb{M}$, then $f\in \mathbb{M}$ (see, e.g., \cite{Z1}). 

For any object $C$ of $\mathcal{C}$, the equivalence relation on the class of $\mathbb{E}$-morphisms with a domain $C$ arises in a natural way. The category $\mathcal{C}$ is said to be $\mathbb{E}$-co-well-powered if the quotient of the latter class by this equivalence relation is a set \cite{B}. 
\begin{theo}
Let $\mathcal{C}$ be a complete category, and let $(\mathbb{E},\mathbb{M})$ be a factorization system on $\mathcal{C}$ with $\mathbb{E}\subseteq Epi$. Let $\mathcal{C}$ be $\mathbb{E}$-co-well-powered, and let $\mathcal{X}$ be a full replete subcategory of $\mathcal{C}$. The following conditions are equivalent:\vskip+2mm

(i) $\mathcal{X}$ is closed under products and $\mathbb{M}$-subobjects (the latter means that if a morphism $m:X\rightarrow X'$ lies in $\mathbb{M}$ and $X'$ is an object of $\mathcal{X}$, then $X$ also is an object of $\mathcal{X}$);\vskip+2mm 

(ii)  $\mathcal{X}$ is a reflective subcategory of $\mathcal{C}$, and the reflection unit components lie in $\mathbb{E}$;\vskip+2mm

(iii) $\mathcal{X}$ is a reflective subcategory of $\mathcal{C}$ and is closed under $\mathbb{M}$-subobjects.
\end{theo}

\begin{proof} (i)$\Rightarrow$(ii): Let $C$ be an object of $\mathcal{C}$. For any equivalence class $\varpi$ of $\mathbb{E}$-morphisms with domain $C$ and codomains in $\mathcal{X}$, consider a representative $e_{\varpi}$. Consider the product $X=\prod_{\varpi} codom \ e_{\varpi}$. Note that $X$ is an object of the subcategory $\mathcal{X}$. The family of morphisms $e_{\varpi}$ induces a morphism $\varrho:C\rightarrow X$. Consider the $(\mathbb{E},\mathbb{M})$-factorization $me$ of $\varrho$. Obviously the codomain $X'$ of the morphism $e$ is an object of $\mathcal{X}$. We assert that $e$ is a universal morphism from the object $C$ to the subcategory $\mathcal{X}$. Indeed, consider any $f:C\rightarrow Y$ with $Y$ in $\mathcal{X}$. Let $m'e'$ be the  $(\mathbb{E},\mathbb{M})$-factorization of $f$ with  $e':C\twoheadrightarrow Y'$. Then $Y'$ also is an object of $\mathcal{X}$. Without loss of generality we can assume that the representative $e_{\varpi}$ of $e'$'s equivalence class ${\varpi}$ is $e'$ itself. Therefore $\pi_{e'}\varrho=e'$, where $\pi_{e'}$ is the canonical projection $X\rightarrow codom \ e'$. Then $(m'\pi_{e'} m)e=f$.

(ii)$\Rightarrow$(iii): Consider a morphism $m:X\rightarrow X'$ from $\mathbb{M}$ with $X'$ in the subcategory $\mathcal{X}$. Let $\eta$ be the unit of the reflection $r$. There is a morphism $h:r(X)\rightarrow X'$ such that $m=h\eta_X$. Then $\eta_X\in \mathbb{M}$. This implies that $\eta_X$ is an isomorphism, and hence $X$ is an object of $\mathcal{X}$.

The implication (iii)$\Rightarrow$(i) is well-known.
\end{proof}

Note that Theorem 4.1 implies `old theorem' since any complete well-powered category has the factorization system $(Strong Epi, Mono)$ (see, e.g., \cite{B}).

\section{Left hereditary left perfect right coherent rings}

We begin this section with the following well-known theorem.

\begin{theo} (Chase \cite{C})
A ring is left perfect and right coherent if and only if the direct product of an arbitrary small family of projective modules is projective.
\end{theo} 

Recall that $\Lambda$-$Mod_{St}$ denotes the full subcategory of stable modules and  $\Lambda$-$Mod_{Pr}$ denotes the full subcategory of projective modules of $\Lambda$-$Mod$.

\begin{theo}
For a ring $\Lambda$, the conditions (i)-(xv) below are equivalent:

\vskip+2mm

(i) the ring $\Lambda$ is left hereditary, left perfect and right coherent;\vskip+2mm

(ii) the subcategory $\Lambda$-$Mod_{Pr}$ of $\Lambda$-$Mod$ is epi-reflective;
\vskip+2mm 

(iii) the ring $\Lambda$ is left hereditary and the subcategory $\Lambda$-$Mod_{Pr}$ of $\Lambda$-$Mod$ is reflective;
\vskip+2mm 

(iv) the ring $\Lambda$ is left hereditary and, moreover, for any module $M$, there is a representation $M\approx K\oplus P$ such that $P$ is projective and if $M\approx K'\oplus P'$ is a representation with projective $P'$, then there is a projective module $Q$ such that $P\approx P'\oplus Q$, $K'\approx K\oplus Q$, and the following equalities hold:
$$i_Pi_Q=i_{K'}i'_{Q},\; \pi_Pi_{K'}=i_Q\pi'_{Q},\; \pi'_{P'} \pi_P=\pi_{P'}, \; i_{K'}i'_{K}=i_K;$$
\noindent here the letter $i$ with an index (and possibly with the prime) denotes the canonical embedding of the module indicated in the index; similarly the letter $\pi$ with an index (and possibly with the prime) denotes the relevant canonical projection:
\begin{equation}
\begin{xymatrix}
{K'\ar@<-0.5ex>@{->>}[rr]_{\pi'_Q}\ar@{ >->}[dr]^{i_{K'}}&&Q\ar@{ >->}[d]^{i_Q}\ar@<-0.5ex>@{ >->}[ll]_{i'_Q}\\
K\ar@{ >->}[u]^{i'_K}\ar@<-0.5ex>@{ >->}[r]_{i_K}&M\ar@<-0.5ex>@{->>}[l]_{\pi_K}\ar@<-0.5ex>@{->>}[r]_{\pi
_P}\ar@{->>}[dr]_{\pi_{P'}}&P\ar@{->>}[d]^{\pi'_{P'}}\ar@<-0.5ex>@{ >->}[l]_{i_P}\\
&&P'}
\end{xymatrix}
\end{equation}
\vskip+2mm

(v) the ring $\Lambda$ is left hereditary and, moreover, for any module $M$, there is a representation $M\approx K\oplus P$ such that $P$ is projective and if $M\approx K'\oplus P'$ is a representation with projective $P'$, then there is a monomorphism $m:K\rightarrowtail K'$ with $i_{K'}m=i_K$;\vskip+2mm

(vi) the ring is left hereditary, and any module can be represented as the direct sum of a stable module and a projective module;\vskip+2mm

\vskip+2mm
(vii) the ring is left hereditary, and any module can be represented as the direct sum of a stable module and a projective module; such a representation is unique 
in the following sense: if $M\approx S\oplus P$ and $M\approx S'\oplus P'$ with stable $S$ and $S'$, and projective $P$ and $P'$, then there is an isomorphism $S\rightarrow S'$ that is compatible with the canonical embeddings, and an isomorphism $P\rightarrow P'$ that is compatible with the canonical projections;
\vskip+2mm

(viii) the ring is left hereditary, and, moreover, any module can be represented as the direct sum of a stable module and projective indecomposable modules;\vskip+2mm

(ix) the ring is left hereditary, and, moreover, any module $M$ can be represented as the direct sum of a stable module $S$ and projective indecomposable modules $P_i$ $(i\in \mathfrak{A})$; such a representation is unique in the following sense: if $M$ is the direct sum of a stable module $S'$ and indecomposable projective modules $P'_i$ $(i\in \mathfrak{B})$, then there is an isomorphism $S\rightarrow S'$ that is compatible with the canonical embeddings, a bijection $\sigma:\mathfrak{A}\rightarrow \mathfrak{B}$ and isomorphisms $P_i\rightarrow P'_{\sigma(i)}$ ($i\in \mathfrak{A}$) such that the induced isomorphism $\oplus_{i\in \mathfrak{A}} P_i\rightarrow \oplus_{i\in \mathfrak{B}} P'_i$ is compatible with the canonical projections $M\twoheadrightarrow \oplus_{i\in \mathfrak{A}} P_i$ and $M\twoheadrightarrow \oplus_{i\in \mathfrak{B}} P'_i$; moreover, the decomposition $\oplus_{i\in \mathfrak{B}} P_i$ complements direct summands of the module $N= \oplus_{i\in \mathfrak{A}} P_i$;\vskip+2mm

(x) the ring $\Lambda$ is left hereditary and, moreover, any module $M$ has a representation $M= \textbf{R}(M)\oplus P$ with a projective module $P$;\vskip+2mm

(xi) the ring $\Lambda$ is left hereditary, and, moreover, for any module $M$, the module $M/\textbf{R}(M)$ is projective;\vskip+2mm

(xii) the ring $\Lambda$ is left hereditary, and, moreover, for any module $M$, $\textbf{R}(M)=0$ if and only if $M$ is projective;\vskip+2mm

(xiii) the ring is left hereditary, and, moreover, the class $(\Lambda$-$Mod_{St})^{\perp}$ coincides with that of all projective modules;\vskip+2mm

(xiv) the pair of module classes $(\Lambda$-$Mod_{St},\Lambda$-$Mod_{Pr})$ is a pre-torsion theory;\vskip+2mm

(xv) $\Lambda$-$Mod_{St}=(\Lambda$-$Mod_{Pr})^{\top}$ and $\Lambda$-$Mod_{Pr}=(\Lambda\text{-}Mod_{St})^{\bot}$. \vskip+2mm

\end{theo}
\begin{proof}

The equivalences (i)$\Leftrightarrow$(ii)$\Leftrightarrow$(iii) immediately follow from Theorem 5.1 and Theorem 4.1 (take $\mathbb{E}=Epi$ and $\mathbb{M}=Mono$). \vskip+1mm

(ii)$\Rightarrow$(iv): First note that the ring $\Lambda$ is left hereditary. Let now \begin{equation}
\textbf{r}:\Lambda\text{-}Mod\rightarrow \Lambda\text{-}Mod_{Pr}
\end{equation}
 \noindent be the reflector, and $\eta$ be the unit of the reflection. Let $M\approx K\oplus \textbf{r}(M)$ be the representation induced by the homomorphism $\eta_M:M\twoheadrightarrow \textbf{r}(M)$. Consider an arbitrary representation $M\approx K'\oplus P'$ with projective $P'$. There is a homomorphism $\varphi:\textbf{r}(M)\rightarrow P'$ with $\varphi\eta_M=\pi_{P'}$ (see diagram (5.3) below). This implies that there is $m:K\rightarrow K'$ with $i_{K'}m=i_K$. Moreover, $\varphi$ is an epimorphism, and hence $\textbf{r}(M)=P'\oplus Q$ for the submodule $Q=Ker\ \varphi$ of $\textbf{r}(M)$.  This, in particular, implies that $Q$ is projective. 

Further, we have $\varphi=\varphi\eta_Mi_{\textbf{r}(M)}=\pi_{P'}i_{\textbf{r}(M)}$. Therefore there is a homomorphism $u:Q\rightarrow K'$ with $i_{K'}u=i_{\textbf{r}(M)}i_Q$; and since $\varphi\eta_Mi_{K'}=0$, there is a homomorphism $v:K'\rightarrow Q$ such that $i_{Q}v=\eta_Mi_{K'}$. Let $e=\pi_Ki_{K'}$:
\begin{equation}
\begin{xymatrix}
{K'\ar@<-0.5ex>@{-->>}[rr]_{v}\ar@{ >->}[dr]^{i_{K'}}\ar@<0.5ex>@{->>}[d]^{e}&&Q\ar@{ >->}[d]^{i_Q}\ar@<-0.5ex>@{ >-->}[ll]_{u}\\
K\ar@<0.5ex>@{ >-->}[u]^{m}\ar@<-0.5ex>@{ >->}[r]_{i_K}&M\ar@<-0.5ex>@{->>}[l]_{\pi_K}\ar@<-0.5ex>@{->>}[r]_{\eta_M}\ar@{->>}[dr]_{\pi_{P'}}&\textbf{r}(M)\ar@{-->>}[d]^{\varphi}\ar@<-0.5ex>@{ >->}[l]_{i_{\textbf{r}(M)}}\\
&&P'}
\end{xymatrix}
\end{equation}
 Then the quantiple $(K', K, Q, m, e, u, v)$ is a biproduct. Indeed, the equalities $em=1_K$ and $eu=0$ are obvious. Moreover, we have $i_Qvu=i_Q$. This implies that $vu=1_Q$. Similarly, $i_{Q} vm=0$, whence $vm=0$. Finally, it is easy to verify that $i_{K'} (me+uv)=(i_K\pi_K+i_{\textbf{r}(M)}\eta_M)i_{K'}=i_{K'}$, and hence $me+uv=1_{K'}$. Thus $K'\approx K\oplus Q$.\vskip+1mm
  
The implication (iv)$\Rightarrow$(v) is obvious.\vskip+1mm
 
(v)$\Rightarrow$(vi): The module $K$ is obviously stable.\vskip+1mm

(vi)$\Rightarrow$(ii): Let $M=S\oplus P$ with a stable module $S$ and a projective module $P$. The canonical projection $\pi_P:S\oplus P\twoheadrightarrow P$ is a universal morphism from $M$ to the subcategory $\Lambda$-$Mod_{Pr}$. Indeed, if $f:M\rightarrow P'$ is a homomorphism with a projective module $P'$, then, according to Lemma 3.2, we have $fi_S=0$. Thus there is a homomorphism $\psi:P\rightarrow P'$ with $\psi \pi_P=f$.\vskip+1mm

For the implication (iv)$\Rightarrow$(vii) it suffices to observe that the module $K$ in the condition (iv) is stable. The implication (vii)$\Rightarrow$(vi) is obvious. 

The implication (i)$\Rightarrow$(viii) follows from Theorem 6 of \cite{AF} (and the fact that (i) implies (vi), shown above). The implication (viii)$\Rightarrow$(vi) is obvious.\vskip+1mm

The implication (i)$\Rightarrow$(ix) follows from Theorems 2, 6 of \cite{AF} (and the fact that (i) implies (vii), as shown above). The implication (ix)$\Rightarrow$(vi) is obvious.\vskip+1mm

(iv)$\Rightarrow$(x): We already know that the module $K$ in the representation $M\approx K\oplus P$ is stable. Assume that $S$ is a stable submodule of $M$ with the embedding $i:S\rightarrowtail M$.  By Lemma 3.2, $\pi_P i=0$, and hence there is a monomorphism $\varphi:S\rightarrowtail K$ such that $i_K\varphi =i$. Thus $K$ is the largest stable submodule of $M$.\vskip+1mm

The implications (x)$\Rightarrow$(xi)$\Rightarrow$(xii) are obvious. 

The equivalence (xii)$\Leftrightarrow$(xiii) follows from Proposition 3.4(e), while the implication (xiii)$\Rightarrow$(xiv) follows from Proposition 3.4(b). \vskip+1mm

The implication (xiv)$\Rightarrow$(iii) follows from the well-known facts from torsion theory. They imply also the equivalence (xiv)$\Leftrightarrow$(xv).
\end{proof}

\begin{cor}
Let $\Lambda$ be a left hereditary left perfect and right coherent ring. Then, for any module $M$, the submodule $\textbf{R}(M)$ is smallest among submodules $N$ with projective $M/N$.

\begin{rem} Let $\Lambda$ be a left hereditary left perfect and right coherent ring. The stable-projective representation of a module $M$, in fact, has the form
 \begin{equation}
M\approx \textbf{R}(M)\oplus \textbf{r}(M),
\end{equation} \noindent and hence is functorial (here $\textbf{r}$ is the reflector (5.2)). However, this representation is natural if and only if the ring is semisimple. Indeed, let the diagram
\begin{equation}
\xymatrix{M\ar[r]^-{\approx}\ar[d]_{f}&\textbf{R}(M)\oplus \textbf{r}(M)\ar[d]^{\textbf{R}(f)\oplus \textbf{r}(f)}\\
M'\ar[r]^-{\approx}&\textbf{R}(M')\oplus \textbf{r}(M')}
\end{equation}
\vskip+1mm
\noindent be commutative, for any $f$, and assume that there is a nonprojective module $M'$.  Take an epimorphism $f:M\twoheadrightarrow M'$ with a projective module $M$. The module $\textbf{R}(M)$ is zero. By Theorem 5.2,  $\textbf{R}(M')$ is not zero. Then $\textbf{R}(f)\oplus \textbf{r}(f)$ can not be an epimorphism, and we arrive to the contradiction.\vskip+1mm
\vskip+2mm
\end{rem}

\end{cor}


\begin{theo}
Let $\Lambda$ be a principal ideal domain that is not a field. Then the fundamental theorem on finitely generated $\Lambda$-modules can not be generalized to the case of all $\Lambda$-modules.
\end{theo}

\begin{proof}
Recall that any perfect Noetherian ring is Artinian. Moreover, an Artinian integral domain is a field. Now it suffices to observe that the representation of modules mentioned in the fundamental theorem induces the stable-projective representation as any torsion module over an integral domain is stable. 
\end{proof}

\begin{rem}
In the case of the ring of integers, the statement in Theorem 5.5 is well-known. It follows, for instance, from the Kulikov's criterion for the existence of representations of Abelian $p$-groups in the form of the direct sums of cyclic groups \cite{Ku}. 
\end{rem}


Further, Theorem 5.2 and the main result of the paper \cite{WMJ} (see Theorem on page 8 and Remarks 1, 2) imply
\begin{theo}
For a ring, $\Lambda$, the pair $(\Lambda$-$Mod_{St},\Lambda$-$Mod_{Pr})$ of module classes is a torsion theory if and only if $\Lambda$ is left hereditary and the injective envelope of the module $_\Lambda \Lambda$ is projective.
\end{theo}
\begin{proof}
It suffices to observe that if the ring $\Lambda$ is left hereditary and the injective envelope of the module $_\Lambda \Lambda$ is projective, then $\Lambda$ is left perfect and right coherent \cite{CR} (see also \cite{MZ}).
\end{proof}

Theorem 5.7 immediately implies Corollary 5.8 below. In particular, one obtains the new proof of the equivalence of its conditions (i) and (iii) that is the result by Colby and Rutter \cite{CR}. Note that they also gave the intrinsic description of left hereditary rings satisfying these conditions in the same paper. In \cite{MZ}, Martsinkovsky and the present author gave the characterization of such rings in terms of the projectively stable category of the ring. 

\begin{cor}\label{st2}
Let $\Lambda$ be left hereditary. Then the following conditions are equivalent:
\vskip+1mm
(i) the injective envelope of the module $_\Lambda \Lambda$ is projective;
\vskip+1mm
(ii) the subcategory $\Lambda$-$Mod_{Pr}$ is closed under essential extensions;
\vskip+1mm
(iii) the subcategory $\Lambda$-$Mod_{Pr}$ is closed under injective envelopes.
\end{cor}\vskip+2mm

\begin{rem}
Let $\Lambda$ be a left hereditary ring and the injective envelope of the module $_\Lambda \Lambda$ be projective. It is easy to observe that torsion theory $(\Lambda$-$Mod_{St},\Lambda$-$Mod_{Pr})$ coincides with the largest torsion theory for which the module $_\Lambda \Lambda$ is torsionfree, considered by Utumi in \cite{U} and Lambek in \cite{L2}. The modules which are divisible with respect to this torsion theory are precisely the injective modules. The divisible envelope of any projective module is its injective envelope.
\end{rem}

Further implications of Theorem 5.7 will be given in our next paper.

\section{Left hereditary left Noetherian rings}

The categorical approach employed in Section 5 for characterizing left hereditary left perfect right coherent rings leads to a new proof of the He's characterization of rings over which all modules have unique injective-costable decompositions. Before we give the proof, recall the following well-known theorem.


\begin{theo} (Bass \cite{Ba1})
A ring is left Noetherian if and only if the direct sum of an arbitrary small family of injective modules is injective.
\end{theo}
Recall that  $\Lambda$-$Mod_{Inj}$ denotes the full subcategory of injective modules and $\Lambda$-$Mod_{Cost}$ denotes the full subcategory of costable modules of $\Lambda$-$Mod$. 
\begin{theo}
For a ring $\Lambda$, the conditions (i)-(xv) below are equivalent:

\vskip+2mm

(i) the ring $\Lambda$ is left hereditary and left Noetherian;\vskip+2mm

(ii) the subcategory $\Lambda$-$Mod_{Inj}$ of $\Lambda$-$Mod$ is mono-coreflective;
\vskip+2mm 

(iii) the ring $\Lambda$ is left hereditary and the subcategory $\Lambda$-$Mod_{Inj}$ of $\Lambda$-$Mod$ is coreflective;
\vskip+2mm 

(iv) the ring $\Lambda$ is left hereditary and, moreover, for any module $M$, there is a representation $M\approx I\oplus K$ such that $I$ is injective and if $M\approx I'\oplus K'$ is a representation with injective $I'$, then there is an injective module $J$ with $I\approx I'\oplus J$, $K'\approx K\oplus J$, and such that the following equalities hold:
$$\pi'_J\pi_{K'}=\pi_J\pi_I,\; \pi_{K'} i_I=i'_J\pi_J,\; i_Ii'_{I'}=i_{I'}, \pi'_K\pi_{K'}=\pi_{K};$$
\noindent here the letter $i$ with an index (and possibly with prime) denotes the relevant canonical embedding, while the letter $\pi$ with an index (and possibly with prime) denotes the relevant canonical projection:
\begin{equation}
\begin{xymatrix}
{I'\ar@{>->}[d]_{i'_{I'}}\ar@{>->}[dr]^{i_{I'}}&&\\
I\ar@<-.5ex>@{ >->}[r]_{i_I}\ar@{->>}[d]_{\pi_J}&M\ar@<-.5ex>@{->>}[l]_{\pi_I}\ar@<.5ex>@{->>}[r]^{\pi_K}\ar@{->>}[dr]_{\pi_{K'}}&K\ar@<.5ex>@{ >->}[l]^{i_K}\\
J\ar@<-.5ex>@{ >->}[rr]_{i'_J}&&K'\ar@<-.51ex>@{->>}[ll]_{\pi'_J}\ar@<.5ex>@{->>}[u]_{\pi'_K}}
\end{xymatrix}
\end{equation}
\vskip+2mm

(v) the ring $\Lambda$ is left hereditary and, moreover, for any module $M$, there is a representation $M\approx I\oplus K$ such that $I$ is injective and if $M\approx I'\oplus K'$ is a representation with injective $I'$, then there is an epimorphism $e:K'\twoheadrightarrow K$ that is compatible with the canonical projections;\vskip+2mm

(vi) the ring is left hereditary, and any module can be represented as the direct sum of an injective module and a costable module;\vskip+2mm

(vii) the ring is left hereditary, and any module can be represented as the direct sum of an injective module and a costable module; such a representation is unique in the following sense: if $M\approx I\oplus C$ and $M\approx I'\oplus C'$ with injective $I$ and $I'$ and costable $C$ and $C'$, then there is an isomorphism $I\rightarrow I'$ that is compatible with the canonical embeddings an isomorphism $C\rightarrow C'$ that is compatible with the canonical projections;\vskip+2mm

(viii) the ring is left hereditary, and, moreover, any module can be represented as the direct sum of injective indecomposable modules and a costable module;\vskip+2mm

(ix) the ring is left hereditary, and, moreover, any module can be represented as the direct sum of injective indecomposable modules $I_i$ ($i\in \mathfrak{A}$) and a costable module $C$; such a representation is unique in the following sense: if $M$ is the direct sum of injective indecomposable modules $I'_i$ ($i\in \mathfrak{B}$) and a costable module $C'$, then there is an isomorphism $C\rightarrow C'$ that is compatible with the canonical projections, a bijection $\sigma:\mathfrak{A}\rightarrow \mathfrak{B}$ and isomorphisms $I_i\rightarrow I'_{\sigma(i)}$ ($i\in \mathfrak{A}$) such that the induced isomorphism $\oplus_{i\in \mathfrak{A}} I_i\rightarrow \oplus_{i\in \mathfrak{B}} I'_i$ is compatible with the canonical embeddings $\oplus_{i\in \mathfrak{A}} I_i\rightarrowtail M$ and $\oplus_{i\in \mathfrak{B}} I'_i\rightarrowtail M$; moreover, the decomposition $\oplus_{i\in \mathfrak{B}} I_i$ complements direct summands of the module $N= \oplus_{i\in \mathfrak{A}} I_i$;\vskip+2mm

(x) the ring $\Lambda$ is left hereditary and, moreover, any module $M$ has a representation $M= I\oplus M/\textbf{R}'(M)$ with an injective module $I$;\vskip+2mm

(xi) the ring $\Lambda$ is left hereditary, and, moreover, for any module $M$, the module $\textbf{R}'(M)$ is injective;\vskip+2mm

(xii) the ring $\Lambda$ is left hereditary, and, moreover, for any module $M$, $\textbf{R}'(M)=M$ if and only if $M$ is injective;\vskip+2mm

(xiii) the ring is left hereditary, and the class $(\Lambda$-$Mod_{Cost})^{\top}$ coincides with that of all injective modules;\vskip+2mm

(xiv) the pair of module classes $(\Lambda$-$Mod_{Inj},\Lambda$-$Mod_{Cost})$ is a pre-torsion theory;\vskip+2mm

(xv) $\Lambda$-$Mod_{Cost}=(\Lambda$-$Mod_{Inj})^{\bot}$ and $\Lambda$-$Mod_{Inj}=(\Lambda\text{-}Mod_{Cost})^{\top}$. \vskip+2mm
\end{theo}

\begin{proof} The equivalence (i)$\Leftrightarrow$(ii)$\Leftrightarrow$(iii) immediately follows from Theorem 6.1 and the dual of Theorem 4.1. The implication (i)$\Rightarrow$(viii) follows from Theorem 2.5 by Matlis \cite{Ma}, and the implication (vii)$\Rightarrow$(ix) follows from Theorems 2, 8 by Anderson and Fuller \cite{AF}. The proofs of all other implications are similar to those of the corresponding implications in Theorem 5.2.\vskip+1mm
\end{proof}

\begin{cor}
Let $\Lambda$ be a left hereditary and left Noetherian ring. Then the class of i-reduced modules is closed under direct products and extensions.
\end{cor}

\begin{cor}
Let $\Lambda$ be a left hereditary and left Noetherian ring. For any module $M$, the submodule $\textbf{R}'(M)$ is the largest injective submodule of $M$.
\end{cor}

\begin{rem} Let $\Lambda$ be a left hereditary left Noetherian ring. The injective-costable representation of a module $M$, in fact, has the form
 \begin{equation}
M\approx \textbf{c}(M)\oplus (M/\textbf{R}'(M)),
\end{equation} \noindent and hence is functorial. Here $\textbf{c}$ is the coreflector 
$$\Lambda\text{-}Mod\rightarrow \Lambda\text{-}Mod_{Inj}.$$
However, this representation is natural if and only if the ring is semisimple. This can be easily shown similarly to the corresponding statement in Remark 5.3 (for $f$, take the essential extension $M\rightarrowtail I(M)$ of a module $M$).
\end{rem}

Note that the pair $(\Lambda$-$Mod_{Inj},\Lambda$-$Mod_{Cost})$ (appearing in the condition (xiv) of Theorem 6.2) is a torsion theory if and only if the ring $\Lambda$ is semisimple.\vskip+2mm

Finally note that Lemma 3.1, Theorem 5.2, and Theorem 6.2 immediately imply the following
\begin{theo}
If any module over a left hereditary ring has a finite indecomposable decomposition, then $\Lambda$ is left perfect and left Noetherian, and hence it is left Artinian.
\end{theo}

\vskip+2mm
\textit{Author's address:
Dali Zangurashvili, A. Razmadze Mathematical Institute of Tbilisi State University},
\textit{1 Aleksidze Str., Tbilisi 0193, Georgia, e-mail: dali.zangurashvili@tsu.ge}

\end{document}